\documentclass[12pt]{article}
\pdfoutput=1

\newlength{\abstractwidth}
\usepackage{epsf}
\usepackage{color}
\usepackage{graphicx}

\usepackage{amsmath}
\usepackage{amssymb}
\usepackage{latexsym}

\usepackage{tikz, pgf}
\usepackage{tkz-fct}
\usepackage{pgfplots}
\usetikzlibrary{shapes.misc}
\usetikzlibrary{shapes,snakes}
\usetikzlibrary{decorations.pathmorphing}	
\usetikzlibrary{decorations.markings}

\renewcommand{\thefootnote}{\fnsymbol{footnote}}
\renewcommand{\thanks}[1]{\footnote{#1}}
\newcommand{\starttext}{
\setcounter{footnote}{0}
\renewcommand{\thefootnote}{\arabic{footnote}}}

\newcommand{\bea}{\begin{eqnarray}}
\newcommand{\eea}{\end{eqnarray}}
\newcommand{\be}{\begin{eqnarray}}
\newcommand{\ee}{\end{eqnarray}}
\newcommand{\<}{\langle}
\renewcommand{\>}{\rangle}

\def\cE{{\cal E}}

\def\tr{{\rm tr}}

\def\det{{\rm det \,}}

\def\p{\partial}
\def\al{\alpha}
\def\e{\epsilon}

\def\om{\omega}

\def\m{\mu}

\def\sq{{{\sqrt{{\scalebox{0.75}[1.0]{\( - 1\)}}}}}\hskip .01in}
\def\ddb{\partial\bar\partial}

\def\o{\omega}
\def\oKE{\omega_{_{\rm KE}}}
\def\ddbar{\partial\bar\partial}
\def\v{\vskip .1in}

\def\sub{\subseteq}
\def\PSH{{\rm PSH}}
\def\I{\int}
\definecolor{Cyan}{cmyk}{1.,0,0,0}
\definecolor{Magenta}{cmyk}{0,1.,0,0}
\definecolor{Yellow}{cmyk}{0,0,1.,0}
\definecolor{White}{cmyk}{0,0,0,0}
\definecolor{Orange}{cmyk}{0,0.61,0.87,0}
\definecolor{RedOrange}{cmyk}{0,0.77,0.87,0}
\definecolor{Red}{cmyk}{0,1.,1.,0}
\definecolor{Purple}{cmyk}{0.45,0.86,0,0}
\definecolor{Violet}{cmyk}{0.79,0.88,0,0}
\definecolor{Blue}{cmyk}{1,0.5,0,0}
\definecolor{ProcessBlue}{cmyk}{0.96,0,0,0}
\definecolor{GreenYellow}{cmyk}{0.6,0,1.,0}
\definecolor{Black}{cmyk}{0,0,0,1}

\newcommand{\ric}{\mathrm{Ric}}

\newtheorem{theorem}{Theorem}
\newtheorem{lemma}{Lemma}

\begin{document}
\starttext
\setcounter{footnote}{0}

\begin{center}

{\Large \bf K\"ahler-Einstein Metrics and Eigenvalue Gaps
\footnote{Work supported in part by the National Science Foundation under grants DMS-1855947 and DMS-1945869.}}

\vskip 0.2in

{\large Bin Guo, Duong H. Phong, and Jacob Sturm}

\vskip 0.3in

\begin{abstract}

The existence of K\"ahler-Einstein metrics on a Fano manifold is characterized in terms of a uniform gap between $0$ and the first positive eigenvalue of the Cauchy-Riemann operator on smooth vector fields. It is also characterized by a similar gap between 0 and the first positive eigenvalue for Hamiltonian vector fields. The underlying tool is a compactness criteria for suitably bounded subsets of the space of K\"ahler potentials which implies a positive gap.

\end{abstract}
\end{center}

\setcounter{tocdepth}{2} 

\baselineskip=15pt
\setcounter{equation}{0}
\setcounter{footnote}{0}

\section{Introduction}
\setcounter{equation}{0}
\label{sec:1}

Starting with the works of Calabi \cite{Ca} and Yau \cite{Y}, a central problem in K\"ahler geometry has been determining when a complex manifold admits a constant scalar curvature K\"ahler metric in a given K\"ahler class. One of the first obstructions to the existence of a cscK K\"ahler metric is the vanishing of the Futaki invariant, which is a character defined on the Lie algebra of holomorphic vector fields. The Yau-Tian-Donaldson conjecture \cite{Y1,T,D02} (see also \cite{PS} for a review) asserts that the existence of a K\"ahler metric with constant scalar curvature should be equivalent to the algebro-geometric notion of K-stability.
Two recent major advances on this conjecture have been the solution of X.X. Chen, S. Donaldson, and S. Sun \cite{CDS1, CDS2, CDS3} of the case of K\"ahler-Einstein metrics on Fano manifolds, and the more recent works by X.X. Chen and J.R. Cheng \cite{CC1, CC2, CC3} which established the equivalence between the existence of a K\"ahler metric with constant scalar curvature and an analytic notion of $K$-stability.

\smallskip
The K-stability condition of a K\"ahler class is the requirement that the generalized Futaki invariant attached to a test configuration be non-negative, and vanish only if the test configuration is a product. It is only one possible characterization of the existence of a canonical metric, and for both geometric and analytic reasons, it may be useful to have other characterizations as well. In the case of K\"ahler-Einstein metrics, which are the focus of the present paper, a notion of $\delta$-invariant has been proposed by Fujita-Odaka \cite{FO} and Blum and Jonsson \cite{BJ}, and it has been shown by R. Berman, S. Boucksom, and M. Jonsson \cite{BBJ}
that the existence of a K\"ahler-Einstein metric is equivalent to the $\delta$-invariant being greater or equal to $1$. In a different and even earlier direction, it had been shown in \cite{PS, PSSW1, ZZ} that the K\"ahler-Ricci flow converges if the lowest strictly positive eigenvalue of the $\bar\p$ operator on vector fields remains bounded uniformly away from $0$ along the flow. It was suggested there \cite{PS, PSSW1, PSSW2} that it may be possible to characterize the existence of a K\"ahler-Einstein metric in terms of lower bounds for this eigenvalue, and this is the problem which we solve in the present paper.

\smallskip
More precisely,
let $X$ be a compact K\"ahler manifold with $c_1(X)>0$. Fix a reference metric $\o_0\in c_1(X)$. For any $\omega\in c_1(X)$, let $K_{\o_0}(\omega)$ be the $K$-energy of $\omega$ with respect to the reference metric $\omega_0\in c_1(X)$, and $u_\omega$ be the normalized Ricci potential of $\omega$, as defined in (\ref{eqn:ricci potential}) below. We define $\lambda_\omega$ to be the lowest strictly positive eigenvalue of the $\bar\p$ operator on vector fields, i.e.,
\bea
\label{lambda}
\lambda_{\omega}={\rm inf}_{V\in T^{1,0}(X), V\perp_\o H^0(X,T^{1,0})}
{\|\bar\p V\|_\omega^2\over \|V\|_{\omega}^2}
\eea
where the subindex denotes the $L^2$ norms taken with respect to the metric $\omega$, and $\perp_\o$ indicates the  perpendicularity condition  with respect to $\o$.
Let $R_\omega$ be the scalar curvature of $\om$.
For each $A>0$, we introduce the following subset of the space of K\"ahler metrics in $c_1(X)$,
\bea
c_1(X;A)=\{\omega\in c_1(X); \|u_\omega\|_{C^0}
+\|\nabla_\omega u_\omega\|_{C^0}+\|R_\omega\|_{C^0}\leq A,
\ K_{\o_0}(\omega)\leq A\},
\eea
and the corresponding {\it eigenvalue gap} for the set $c_1(X;A)$ by
\bea
\label{lambda-A}
\lambda(X;A)=\inf_{\omega\in c_1(X;A)}\lambda_\omega
\eea
Then we have the following characterizations of the existence of a K\"ahler-Einstein metric:

\begin{theorem}
\label{main}
Let $X$ be a compact K\"ahler manifold with $c_1(X)>0$ and vanishing Futaki invariant. Then $X$ admits a K\"ahler-Einstein metric if and only if $\lambda(X,A)>0$ for any $A>0$.
\end{theorem}

Note that although the definition of the $K$-energy requires a choice of reference metric $\om_0$, under a change of reference metric, it just shifts by a constant.
Thus the above condition is invariant under a change of reference metric, as it should be.

\medskip

To explain the second characterization, we recall the following observations due by Futaki \cite{Fut} (see also \cite{PSSW1}, Lemma 2). For any metric $\o\in c_1(X)$, the differential operator operating on smooth functions
$L_\o f\ = \ -g^{i\bar j}\partial_i\partial_{\bar j}f + 
g^{i\bar j}\partial_iu_\o\,\partial_{\bar j}f -  f$ is non-negative,
and its kernel is the space of functions $f$ with $\nabla f$ a holomorphic vector field.
Let $\mu_\o$ be the smallest positive eigenvalue of $L_\omega$. Then the corresponding eigenfunctions $f$ satisfy the identity
\bea
\label{futaki}
\int_X |\bar\nabla\bar\nabla f|^2 e^{-u_\o}\omega^n
= \m_\omega \int_X |\bar\nabla f|^2e^{-u_\o}\omega^n.
\eea
Moreover
\bea
\mu_\o={\rm inf}_{f\in C^\infty(X), \,\int_X fe^{-u_\o}\o^n=0}{\int_X|\bar\nabla\bar\nabla f|^2 e^{-u_\o}\omega^n\over
\int_X|\bar\nabla f|^2e^{-u_\o}\omega^n}
\eea
We introduce, in analogy with (\ref{lambda-A}), the eigenvalue gap for Hamiltonian vector fields by
\bea
\label{mu-A}
\mu(X;A)=\inf_{\omega\in c_1(X;A)}\mu_\omega.
\eea

\begin{theorem}
\label{main1}
Let $X$ be a compact K\"ahler manifold with $c_1(X)>0$ and vanishing Futaki invariant. Then $X$ admits a K\"ahler-Einstein metric if and only if $\mu(X,A)>0$ for any $A>0$.
\end{theorem}
\smallskip

For each $A>0$, we have the easy bound $\mu(X;A)\geq c_A\,\lambda(X;A)$ for some positive constant $c_A$. Thus the condition $\lambda(X;A)>0$ in Theorem \ref{main} implies the condition $\mu(X;A)>0$ in Theorem \ref{main1}. However, there does not appear to be a direct way to show that they are equivalent. 

\bigskip

We now describe briefly our approach. One direction in Theorem \ref{main} and Theorem \ref{main1} is known, by combining the work of Perelman on the K\"ahler-Ricci flow with the convergence results of \cite{PSSW1} and \cite{ZZ}. The main problem is to establish the other direction, namely that the existence of a K\"ahler-Einstein metric on $X$ implies that the gaps $\lambda(X;A)$ and $\mu(X;A)$ are strictly positive for any $A>0$.
For each fixed $\om$, the eigenvalues  
$\lambda_\om$  and $\m_\om$ are positive by definition. So the desired statement can be interpreted as a compactness statement with respect to a suitable topology. Our strategy for such a statement is to view the K\"ahler potential $\varphi$ of a metric $\om\in c_1(X,A)$ as the solution of a Monge-Amp\`ere equation with right hand side depending on the Ricci potential $u_\om$.  The $C^\al$ estimates are derived by combining the theorem of Skoda-Zeriahi \cite{Zer} with that of Kolodziej \cite{Kol} following the idea of Guedj \cite{BEG}. Then the
$C^{3,\alpha}$ priori estimates can be obtained by combining methods for the Monge-Amp\`ere equation together with the recent techniques introduced by Chen-Cheng \cite{CC1} for the constant scalar curvature problem. Next, the $C^{2,\alpha}$ bounds imply the uniform equivalence of the metrics. This implies in turn uniform estimates of the corresponding eigenvalues on vector fields, using the arguments of \cite{PSSW0} to handle the orthogonality condition with different metrics to 
holomorphic vector fields. The desired theorems follow.

\section{$C^{1,\alpha}$ estimates for metrics in $c_1(X;A)$}
\setcounter{equation}{0}

First we set up the equation. Let $n$ be the dimension of $X$. If $\om$ is any metric in $c_1(X)$, we define its Ricci potential $u_\om$ by
\bea\label{eqn:ricci potential}
Ric(\om)-\om=-i\p\bar\p u_\omega,
\qquad
\int_X e^{-u_\omega}\om^n=\int_X\om^n
\eea
where $Ric(\om)=-i\p\bar\p\log\om^n$ is its Ricci curvature form. 
Fix now a reference metric $\om_0\in c_1(X)$, and let $Ric(\om_0)$ and $u_0$ are its Ricci form and Ricci potential, respectively. We can then write $\om=\om_0+i\p\bar\p\varphi$,
where $\varphi$ is normalized to satisfy ${\rm sup}_X\varphi=0$. Since
\bea
-i\p\bar\p(u_\om-u_0)+i\p\bar\p\varphi=Ric(\om)-Ric(\om_0)=-i\p\bar\p\log{\om^n\over\om_0^n}
\eea
we find that $\varphi$ satisfies the following complex Monge-Amp\`ere equation
\bea
\label{eqn:MA}
(\om_0+i\p\bar\p\varphi)^n=e^{u_\om-u_0-\varphi+c_{\varphi}}\om_0^n
\eea
where $c_{\varphi}$ is a specific constant, which is determined because $\varphi$, $u_\om$, and $u_0$ have all been normalized. It follows from the normalization of $u_\omega$ and $\varphi$ that $c_\varphi\le 0$.

\medskip

\subsection{The $C^\alpha$ estimates on potential}

The first step is  the following $C^\alpha$ estimate.

\begin{lemma}
\label{calpha}
Assume that $X$ admits a K\"ahler-Einstein metric $\o_{KE}$, which we take as reference metric $\o_0=\o_{KE}$. Then there exists $\al>0$ with the following property.
For any $\o\in c_1(X; A)$, there exists an automorphism $g$ of $X$ such that 
$\|\psi-\sup\psi\|_{C^\al(\o_{KE})}\ \leq \ C(A)$, where $g^*\o=\o_{KE}+i\ddb\psi$.
\end{lemma}

{\it Proof.} 
Since $X$ admits a K\"ahler-Einstein metric, we can apply the Moser-Trudinger inequality. An early form of this inequality was first proved in \cite{T},
a sharp version subsequently in \cite{PSSW} in the case of manifolds without holomorphic vector fields, and the full sharp and general version in \cite{Dar-Rub}.  
In this form, it asserts that for any $\o\in c_1(X)$,
there exists $g\in G$ (here $G$ is the automorphism group of $X$) and $\e>0$ depending on  $X$  such that
\bea
A \ \geq \ K_{\o_{KE}}(\o_\phi)\ \geq \ \e J_{\o_{KE}}(g^*\o_\phi) - {1\over \e}
\eea
where $J_{\o_{KE}}$ is the Aubin-Yau functional with reference metric $\omega_{KE}$. 
Thus if we write $g^*\o_\phi = \o_{KE}+i\ddb\psi$ we have
\bea
\psi-\sup\psi\ \in \ S_{A_3}\ = \ \{\theta\in\cE_1(X,\o_{KE})\sub\PSH(X,\o_{KE})\,:\, \sup\theta=0\ \ {\rm and} \ \ J_{\o_{KE}}(\o_\theta)\leq A_3\}
\nonumber
\eea
where $\PSH(X,\o_{KE})$ is the space of plurisubharmonic functions and $\cE_1(X,\o_0)$ is the space of finite energy potentials.
We now claim:
\v
\begin{enumerate}
\item $S_{A_3}$ is compact with respect to the weak $L^1(\o_0^n)$
topology on $\PSH(X,\o_0)$
\item Every element of $S_{A_3}$ has zero Lelong number at $z$ for all $z\in X$.
\item For every $p\geq 1$ there exists $C(p,\o_0, A)$ such that
\be\label{skoda1} \I _Xe^{-p\theta}\o_0^n\ \leq \ C(p,\oKE, A)\ \ \ 
\hbox{for all $\theta\in S_{A_3}$ }.
\ee
\end{enumerate}

\v

These follow as in \cite{BEG} respectively from Lemma 4.13, Proposition 2.13 and Theorem 4.15 (due to Skoda and Zeriahi) of
\cite{Dar}. Next, applying (\ref{eqn:MA}) with $\varphi=\psi-{\rm sup}_X\psi$, we obtain for $\o\in c_1(X;A)$,
\be\label{MA}
(\oKE+\sq\ddb\psi)^n\ = \ e^{-{(\psi-\sup\psi)}+u_\psi+c_\psi}\oKE^n\ \leq \ C(A)e^{-{(\psi-\sup\psi)}}\omega^n_{KE}
\ee
where we have used the fact that $|u_\psi|\leq A$ and $c_\psi\leq 0$.
Now if we apply (\ref{skoda1}) to (\ref{MA}) we obtain that $\|\psi-\sup\psi\|_{C^\al(\o_{KE})}\leq C(A)$ for some $\alpha = \alpha(n, p)\in (0,1)$ by the theorem of Kolodziej
\cite{Kol}.
Q.E.D.

\subsection{$C^{3,\alpha}$ estimates on potentials}

We return to the study of the equation (\ref{eqn:MA}), for a general compact K\"ahler manifold $X$ and reference metric $\om_0$, not necessarily K\"ahler-Einstein. The goal of the present subsection is to establish the following lemma:

\begin{lemma}
\label{lemma:C3}
Let $\varphi$ be a smooth solution of the Monge-Amp\`ere equation
(\ref{eqn:MA}). Assume that $\|\varphi\|_{C^0}\leq A$, $\|u_\omega\|_{C^0}+\|\nabla_\omega u_\omega\|_{C^0}+\|\Delta_\om u_\omega\|_{C^0}\leq A$.
Then for any $\alpha\in (0,1)$, there exists a constant $C = C(n,A,\om_0,\alpha)>0$
{}  so that
\bea
\|\varphi\|_{C^{3,\alpha}(X,\om_0)}\leq C.
\eea
\end{lemma}

It is convenient to set $F=-\varphi+u_\om-u_0+c_\varphi$, so the equation can be written as
\bea\label{eqn:2.6}
(\omega_0+\ddbar \varphi)^n = e^{F} \omega_0^n, \quad \sup_X \varphi = 0,
\eea
and note that $F$ depends on the K\"ahler potential of $\varphi$. To simplify the notation, we shall denote $u_\om$ by just $u$. Under the assumptions of the lemma, both $\varphi$ and $u$ are bounded, so it follows from the fact that $\om$ and $\om_0$ have the same volume that $|c_\varphi|$ is bounded by a constant $C(A,\om_0)$ as well. Thus we have
\bea
0< \frac{1}{C(A,\om_0)}\le e^F \le C(A,\om_0). 
\eea
We divide the proof of the lemma into the following steps:

\begin{enumerate}
\item Apply Chen-Cheng's argument \cite{CC1}  to show $\Delta_{\omega_0}\varphi$ is in $L^p(X,\omega_0)$ for any $p>0$, hence $\varphi\in C^{1,\beta}(X,\omega_0)$ for any $\beta\in (0,1)$ by elliptic estimates. Here $\Delta_{\om_0}$ is the Laplacian with respect to the reference metric $\om_0$.

\item The H\"older continuity of $\varphi$ and the assumption $\|\nabla u\|^2_{C^0(X,\omega)}\le C(A)$  implies that $u\in C^{0,\alpha'}(X,\omega_0)$ (see Lemma \ref{lemma new 3} below).

\item By a theorem of Li-Li-Zhang \cite{LLZ} (which is an improvement of a result of Yu Wang \cite{W}), we get the $C^{2,\alpha'}(X,\omega_0)$ bound for $\varphi$.

\item After we show $u\in C^{1,\alpha''}(X,\omega_0)$ by elliptic estimates, we get the $C^{3,\alpha}$ estimate  for $\varphi$ by differentiating the Monge-Amp\`ere equation \eqref{eqn:MA}. 

\end{enumerate}

We begin by modifying the arguments in Chen-Cheng \cite{CC1} to derive the following estimates:

\begin{lemma}\label{lemma new 1}
There exists a constant $C=C(A, n, \omega_0)>0$ with 
\bea
\sup_X \|\nabla \varphi\|_{C^0(X,\omega_0)}^2\le C.
\eea

\end{lemma}

{\it Proof.} Denote $\Phi = - F  - \lambda \varphi + \frac 1 2 \varphi^2$ with a constant $\lambda>0$ to be chosen later.  We calculate 
\bea
\label{eqn:new 1}
\Delta_\omega\big( e^{\Phi} (|\nabla \varphi|^2_{\omega_0} +3 )\big) & = (|\nabla \varphi|^2_{\omega_0} +3 ) \Delta_\omega e^{\Phi} + 2 e^\Phi Re \langle \nabla \Phi, \bar \nabla |\nabla \varphi|^2_{\omega_0} \rangle_\omega + e^{\Phi} \Delta_\omega |\nabla \varphi|^2_{\omega_0}.
\eea
We consider the first term in \eqref{eqn:new 1}.
\begin{align*}
\Delta_\omega e^{\Phi} & = e^{\Phi}( \Delta_\omega \Phi + |\nabla \Phi|^2 _\omega )\\
& = e^\Phi \big( (\lambda - 1) \Delta_\omega(-\varphi) - \Delta_\omega u + \Delta_\omega u_0+ \varphi \Delta_\omega \varphi + |\nabla \varphi|^2_\omega + |\nabla \Phi|^2_\omega    \big)\\
& \ge e^\Phi \big( (\lambda - 1 - \varphi - C_0) \tr_\omega\omega_0 - C + |\nabla \varphi|^2_\omega + |\nabla\Phi|^2_\omega    \big)
\end{align*}
where we used the assumption that $|\Delta_\omega u|\le C(A)$, and $C_0 = C_0(\omega_0)>0$ is a constant satisfying $-C_0 \omega_0\le i\ddbar u_0\le C_0\omega_0$.

\smallskip

To deal with the third term in \eqref{eqn:new 1}, we introduce a normal coordinates system for $\omega_0$ at the maximum point $x_0\in X$ of $e^{\Phi}(3 + |\nabla \varphi|^2_{\omega_0})$ such that $g_0=(\tilde g_{i\bar j}) = (\delta_{ij})$ and $dg_0 = 0$ at the point. Moreover, $\omega = (g_{i\bar i}\delta_{ij})$ is diagonal at $x_0$. We calculate at $x_0$,
\begin{align*}
\Delta_\omega |\nabla \varphi|^2 _{\omega_0} & = g^{p\bar p}\frac{\partial^2}{\partial z_p\partial \bar z_p}(\tilde g^{i\bar j} \varphi_{\bar j}\varphi_i)\\
& = g^{p \bar p}\frac{\partial^2 \tilde g^{i\bar j}}{\partial z_p\partial \bar z_p} \varphi_{\bar j}\varphi_i + g^{p\bar p} \varphi_{\bar i p} \varphi_{i \bar p} + g^{p\bar p } \varphi_{\bar i\bar p}\varphi_{i p} + g^{p\bar p} \phi_i \frac{\partial^2 \varphi_{\bar i}}{\partial z_p\partial \bar z_p} + g^{p\bar p} \varphi_{\bar i} \frac{\partial^2 \varphi_{i}}{\partial z_p\partial \bar z_p}\\
& =\tilde R_{j\bar k p\bar p} g^{p\bar p} \varphi_k\varphi_{\bar j} + g^{p\bar p} \varphi_{\bar i p} \varphi_{i \bar p} + g^{p\bar p } \varphi_{\bar i\bar p}\varphi_{i p} + g^{p\bar p} \varphi_i \frac{\partial^2 \varphi_{\bar i}}{\partial z_p\partial \bar z_p} + g^{p\bar p} \varphi_{\bar i} \frac{\partial^2 \varphi_{i}}{\partial z_p\partial \bar z_p}\\
&\ge  - C_1 \tr_\omega \omega_0 |\nabla \varphi|^2_{\omega_0} + g^{p\bar p} \varphi_{\bar i p} \varphi_{i \bar p} + g^{p\bar p } \varphi_{\bar i\bar p}\varphi_{i p}  + 2 Re( \varphi_{i} F_{\bar i}  ),
\end{align*}
where $\tilde R_{i\bar j k\bar l}$ is the bisectional curvature of the metric $g_0$, $-C_1$ is a lower bound of $\tilde R_{i\bar j k\bar l}$,  and in the last inequality we have used the equation below by taking derivatives on both sides of \eqref{eqn:2.6}
$$g^{p\bar p} \frac{\partial^2\varphi_{i}}{\partial z_p\p \bar z_p} = F_i,\,\text{at }x_0.$$
Therefore, we get
\bea
\label{eqn:new 2}
\Delta_\omega\big( e^{\Phi} (|\nabla \varphi|^2_{\omega_0} +3 )\big) 
&\ge&   
e^\Phi \Big\{ (|\nabla \varphi|^2_{\omega_0} +3 )  \Big( (\lambda - 1 - \varphi - C_0) \tr_\omega\omega_0 - C + { |\nabla \varphi|^2_\omega + |\nabla\Phi|^2_\omega }   \Big)\nonumber\\
&&
\quad - C_1 \tr_\omega \omega_0 |\nabla \varphi|^2_{\omega_0} + g^{p\bar p} \varphi_{\bar i p} \varphi_{i \bar p} + {g^{p\bar p } \varphi_{\bar i\bar p}\varphi_{i p}}
\nonumber\\
&&
\quad  + 2 Re( \varphi_{\bar i} F_{i}  ) + 2 Re\big(g^{i\bar i} \Phi_i ( \varphi_j\varphi_{\bar j \bar i} + \varphi_{j\bar i} \varphi_{\bar j}  )    \big)\Big\}.
\eea
The last two terms are equal to  (note that at $x_0$, $\varphi_{j\bar i} = \varphi_{i\bar i} \delta_{ij} = (g_{i\bar i} - 1) \delta_{ij}$)
\begin{align*}
 & 2 Re\big( - \Phi_i \varphi_{\bar i} - (\lambda - \varphi) |\nabla \varphi|^2_{\omega_0} + g^{i\bar i} \Phi_i (g_{i\bar i} -1  )\varphi_{\bar i}  \big)\\ 
=  &  - 2(\lambda - \varphi) |\nabla \varphi|^2_{\omega_0}  - 2 Re \big(    \langle \nabla\Phi, \bar \nabla \varphi\rangle_{\omega}  \big)\\
\ge & -2 (\lambda - \varphi) |\nabla \varphi|^2_{\omega_0} - |\nabla \varphi|^2_{\omega} - |\nabla \Phi|_\omega^2,
\end{align*}
the last two terms on the RHS can be absorbed by the corresponding terms in the first line on the right hand side in \eqref{eqn:new 2}, while
$$
2 Re (g^{i\bar i}\Phi_i \varphi_j \varphi_{\bar i \bar j}) \ge - g^{i\bar i}\Phi_i\Phi_{\bar i} \varphi_{j}\varphi_{\bar j} - g^{i\bar i}\varphi_{i j}\varphi_{\bar i \bar j} = -|\nabla \varphi|^2_{\omega_0} |\nabla \Phi|^2_{\omega} - g^{i \bar i} \varphi_{i j}\varphi_{\bar i\bar j},  
$$ 
and the right hand side above can also be absorbed by terms in the first and second lines of the right hand side in  \eqref{eqn:new 2}. So we get by combining the above that at $x_0$
\begin{equation}\label{eqn:new 3}\begin{split}
0\ge & \Delta_\omega\big( e^{\Phi} (|\nabla \varphi|^2_{\omega_0} +3 )\big)\\ \ge  &   e^\Phi \Big\{ (|\nabla \varphi|^2_{\omega_0} + 3)  \Big( (\lambda - 1 - \varphi - C_0 - C_1) \tr_\omega\omega_0 - C + { \frac 1 2 |\nabla \varphi|^2_\omega}   \Big) - 2 (\lambda - \varphi) |\nabla\varphi|_{\omega_0}^2\Big\}\\
\ge & e^\Phi \Big( |\nabla \varphi|_{\omega_0}^2 (\tr_\omega\omega_0 + \frac 12|\nabla \varphi|_\omega^2) - C |\nabla \varphi|_{\omega_0}^2 - C    \Big)\\
\ge & e^\Phi \Big( c(n, A)|\nabla \varphi|_{\omega_0}^{2(1+\frac 1 n)}  - C |\nabla \varphi|_{\omega_0}^2 - C    \Big)
\end{split}\end{equation}
where we  choose $\lambda = 2+ \| \varphi\|_{L^\infty} + C_0 +C_1$. In the last step we apply the inequality below which follows from Young's inequality (i.e. $a^{\frac 1 n} b^{\frac{n-1}{n}}\le c(n) ( \frac 12 a + b  )$ for some $c(n)>0$)
\begin{align*}
|\nabla \varphi|_{\omega_0}^2 & \le |\nabla \varphi|_{\omega}^2 \tr_{\omega_0} \omega \le |\nabla \varphi|_\omega^2 (\tr_\omega \omega_0)^{{n-1}} \Big(\frac{\omega^n}{\omega_0^n}\Big)\\
& = e^{{F}} \Big( |\nabla \varphi|_\omega^{\frac 2 n} (\tr_\omega\omega_0)^{\frac{n-1}{n}}  \Big)^{n}\le c(n)  e^{{F}} \Big(\frac 1 2 |\nabla \varphi|_\omega^{2} + \tr_\omega\omega_0  \Big)^{n}\\
&\le C(n, A) \Big(\frac 1 2 |\nabla \varphi|_\omega^{2} + \tr_\omega\omega_0  \Big)^{n}.
\end{align*}
From \eqref{eqn:new 3} we conclude that at $x_0$, $|\nabla \varphi|^2_{\omega_0}\le C(n, A)$. Since $x_0$ is a maximum point of $e^\Phi (|\nabla \varphi|_{\omega_0}^2 + 3)$, we see that $\sup_X|\nabla\varphi|_{\omega_0}^2\le C(n, A)$. The lemma is proved.

\bigskip

We next apply the argument in the proof of Theorem 3.1 of Chen-Cheng \cite{CC1}. In our case the functions $F$ and $\varphi$ are bounded so we can simplify the proof a little bit. 
\begin{lemma} \label{prop 1}
For any $p>0$, there exists a constant $C_p=C(n, A, \omega_0, p)>0$ such that 
$$\int_X (\tr_{\omega_0} \omega)^p\omega_0^n \le C_p.$$

\end{lemma}

\medskip
{\it Proof.} We fix a constant $\alpha\ge 1$ which will be determined later. For notational simplicity, we 
write $\Psi = -\alpha F - \lambda \alpha \varphi$, and calculate $\Delta_\omega (e^{\Psi} \tr_{\omega_0} \omega)$,
\bea
\label{eqn:new 4}
\Delta_\omega (e^\Psi \tr_{\omega_0}\omega)
= e^\Psi \tr_{\omega_0}\omega ( \Delta_\omega \Psi + |\nabla \Psi|_\omega^2  ) + 2 e^\Psi Re\langle \nabla \Psi, \bar \nabla \tr_{\omega_0}\omega \rangle_\omega + e^\Psi \Delta_\omega \tr_{\omega_0} \omega.
\eea
We use a normal coordinates system of $\omega_0$, so that $\omega_0 = (\delta_{ij})$, $d g_{0} = 0$ and $g$ (i.e. $\omega$) is diagonal  at a given point. By the standard calculations as in Yau \cite{Y}, the last term in \eqref{eqn:new 4} satisfies
\begin{equation*}\begin{split}
e^\Psi \Delta_\omega \tr_{\omega_0}\omega & \ge e^\Psi\big( - C_2 \tr_\omega \omega_0\, \tr_{\omega_0}\omega + g^{i\bar i} g^{j\bar j} \varphi_{i\bar jk} \varphi_{j\bar i \bar k} + \Delta_{\omega_0} F - R_{\omega_0}   \big)
\end{split}\end{equation*}
where $-C_2$ is a lower bound of the bisectional curvature of $\omega_0$, $\varphi_{i\bar j k}$ denotes the covariant derivative of $\varphi$ under $\nabla_{\omega_0}$ and $R_{\omega_0}$ is the scalar curvature of $\omega_0$. We cannot apply the usual maximum principle here because apriori $\Delta_{\omega_0} F$ is not bounded.

\smallskip

The second term in \eqref{eqn:new 4} satisfies
\begin{align*}
2e^\Psi  Re \langle \nabla \Psi , \bar \nabla \tr_{\omega_0}\omega\rangle_\omega  & \ge - 2 e^\Psi |\nabla \Psi|_\omega | \nabla \tr_{\omega_0} \omega  |_\omega\\
& \ge - e^\Psi \tr_{\omega_0} \omega |\nabla \Psi|^2_\omega - e^\Psi \frac{|\nabla \tr_{\omega_0}\omega| ^2_\omega }{\tr_{\omega_0} \omega}\\
& \ge - e^\Psi \tr_{\omega_0} \omega |\nabla \Psi|^2_\omega - e^\Psi g^{i\bar i} g^{j\bar j} \varphi_{i\bar j k}\varphi_{j\bar i \bar k},
\end{align*}
where in the last step we use the inequality below as in \cite{Y}
\begin{align*}
| \nabla \tr_{\omega_0} \omega  |^2 _\omega & = \sum_i g^{i\bar i} \big|\sum_k \varphi_{k\bar k i}\big|^2 \le \tr_{\omega_0}\omega \sum_i g^{i\bar i}  \sum_j g^{j\bar j} \varphi_{j \bar j i}\varphi_{\bar j j \bar i} \le \tr_{\omega_0}\omega g^{i\bar i} g^{j\bar j}\, \varphi_{i\bar j k}\varphi_{j\bar i \bar k}.
\end{align*}
The first term in \eqref{eqn:new 4} is 
\bea
e^\Psi \tr_{\omega_0}\omega\, \Delta_\omega \Psi 
&= & e^\Psi \tr_{\omega_0}\omega  \big( \alpha \Delta_\omega \varphi - \alpha \Delta_\omega u + \alpha \Delta_\omega u_0 - \lambda \alpha \Delta_\omega \varphi   \big)\nonumber\\
&\ge & e^\Psi \tr_{\omega_0}\omega \big( (\lambda \alpha - \alpha - C_0) \tr_\omega \omega_0 - C(n, A)\alpha   \big),
\eea
where as before $C_0>0$ satisfies $-C_0\omega_0\le i\ddbar u_0\le C_0\omega_0$. Combining the above inequalities we get
\bea
\label{eqn:new 5}
\Delta_\omega ( e^\Psi \tr_{\omega_0} \omega ) & \ge& e^\Psi \Big( (\lambda \alpha - \alpha - C_0- C_2) \tr_{\omega_0} \omega \tr_\omega \omega_0 - C(n, A )\alpha \tr_{\omega_0} \omega
+ \Delta_{\omega_0} F - R_{\omega_0}   \Big)\nonumber\\
& \ge& e^\Psi \Big( \alpha (\tr_{\omega_0} \omega)^{\frac{n}{n-1}} e^{-\frac F{n-1}}  - C(n, A )\alpha \tr_{\omega_0} \omega + \Delta_{\omega_0} F - R_{\omega_0}   \Big)\nonumber\\
&\ge& e^\Psi \Big( c_0 \alpha (\tr_{\omega_0}\omega)^{\frac n{n-1}} + \Delta_{\omega_0} F - C(n, A) \alpha  \Big)
\eea
where we choose $\lambda = C_0+C_2+2$, $c_0=c_0(n,A,\omega_0)>0$ depends on the lower bound of $e^{-\frac F{n-1}}$, and in the last step we apply Young's inequality $\tr_{\omega_0} \omega \le \varepsilon (\tr_{\omega_0}\omega)^{\frac n{n-1}} + C(\varepsilon)$ for a suitable choice of small $\varepsilon>0$. 

We denote $v: = e^\Psi \tr_{\omega_0}\omega >0$ and for any $p\ge 1$ we have by \eqref{eqn:new 5}
\begin{equation}\label{eqn:new 6}\begin{split}
\Delta_\omega v^p & = p v^{p-1}\Delta_\omega v + p(p-1)v^{p-2} |\nabla v|^2_\omega\\
& \ge  p v^{p-1} e^\Psi \Big( c_0 \alpha (\tr_{\omega_0}\omega)^{\frac n{n-1}} + \Delta_{\omega_0} F - C(n, A) \alpha  \Big) + p(p-1) v^{p-3} e^\Psi |\nabla v|_{\omega_0}^2,
\end{split}\end{equation} 
where in the inequality we have applied the observation that $v |\nabla v|^2_\omega = e^\Psi \tr_{\omega_0} \omega |\nabla v|_\omega^2 \ge e^\Psi |\nabla v|_{\omega_0}^2.$ 
Integrating the inequality \eqref{eqn:new 6} over $X$ against the volume form $\omega^n = e^F \omega_0^n$, we obtain
\bea
\label{eqn:new 7}
&&
\int_X  \Big( v^{p-1} e^{\Psi+F} \big( c_0 \alpha (\tr_{\omega_0}\omega)^{\frac n{n-1}} + \Delta_{\omega_0} F \big) + (p-1) v^{p-3} e^{\Psi+F} |\nabla v|_{\omega_0}^2 \Big)\omega_0^n
\nonumber\\
&&
\qquad
\le C(n, A)\alpha  \int_X v^{p-1} e^{\Psi + F }\omega_0^n.
\eea
To deal with the term involving $\Delta_{\omega_0} F$, we will apply the integration by parts.
We calculate
\bea
\label{eqn:new 8}
\int_X v^{p-1} e^{\Psi + F} \Delta_{\omega_0} F \omega_0^n &= &  \int_X v^{p-1} e^{- (\alpha - 1) F - \lambda \alpha \varphi} \Delta_{\omega_0} F \omega_0^n\nonumber\\
 &= &  \int_X\Big( - (p-1) v^{p-2} e^{-(\alpha - 1) F - \lambda\alpha \varphi} \langle \nabla v, \bar \nabla F \rangle_{\omega_0}\nonumber\\
 &&\quad + v^{p-1} e^{-(\alpha - 1) F - \lambda \alpha \varphi} (\alpha - 1) |\nabla F|_{\omega_0}^2 \nonumber\\
& & \quad + v^{p-1} e^{-(\alpha - 1) F - \lambda \alpha \varphi} \lambda \alpha \langle \nabla \varphi,\bar \nabla F \rangle_{\omega_0}\Big)\omega_0^n.
\eea
The second term in the right hand side of \eqref{eqn:new 8} is good. The first term in \eqref{eqn:new 8} satisfies
\begin{align*}
& \int_X - (p-1) v^{p-2} e^{-(\alpha - 1) F - \lambda\alpha \varphi} \langle \nabla v, \bar \nabla F \rangle_{\omega_0} \\ \ge &  -  \int_X (p-1) v^{p-2} e^{\Psi + F} |\nabla v|_{\omega_0} |\nabla F|_{\omega_0}\\
\ge & -\int_X \frac{\alpha - 1}{ 4} v^{p-1} e^{-(\alpha - 1) F - \lambda \alpha \varphi} |\nabla F|^2_{\omega_0}  - \int_X \frac{(p-1)^2}{\alpha - 1} v^{p-3} e^{\Psi + F} |\nabla v|_{\omega_0}^2\\
\ge & -\int_X \frac{\alpha - 1}{ 4} v^{p-1} e^{-(\alpha - 1) F - \lambda \alpha \varphi} |\nabla F|^2_{\omega_0}  - \int_X v^{p-3} e^{\Psi + F} |\nabla v|_{\omega_0}^2
\end{align*}
if we take $\alpha = \alpha(p)\ge p+2$. These negative terms will be cancelled by the positive terms from \eqref{eqn:new 8} and \eqref{eqn:new 7}. Next we look at the third term on the right hand side of \eqref{eqn:new 8}. By Lemma \ref{lemma new 1} we have a bound on $\sup_X |\nabla \varphi|_{\omega_0}$, and thus
\begin{align*}
&~ \int_X v^{p-1} e^{-(\alpha - 1) F - \lambda \alpha \varphi} \lambda \alpha \langle \nabla \varphi,\bar \nabla F \rangle_{\omega_0}\\
\ge & ~  - C \lambda \alpha \int_X v^{p-1} e^{-(\alpha - 1) F - \lambda \alpha \varphi} |\nabla F|_{\omega_0}\\
\ge & ~ -\frac{\alpha - 1}{ 4} \int_X v^{p-1} e^{-(\alpha -1) F - \lambda \alpha \varphi} |\nabla F|_{\omega_0}^2 - \frac{C\alpha^2}{\alpha - 1} \int_X v^{p-1} e^{-(\alpha - 1) F - \lambda\alpha \varphi}\\
\ge & ~ -\frac{\alpha - 1}{ 4} \int_X v^{p-1} e^{-(\alpha -1) F - \lambda \alpha \varphi} |\nabla F|_{\omega_0}^2 - C\alpha \int_X v^{p-1} e^{-(\alpha - 1) F - \lambda\alpha \varphi}  .
\end{align*}
Plugging the above inequalities into \eqref{eqn:new 7} and re-organizing, it follows that 
\begin{align*}
\int_X c_0\alpha v^{p-1}e^{\Psi + F} (\tr_{\omega_0}\omega)^{\frac n{n-1}} \omega_0^n \le C(n, A) \alpha \int_X v^{p-1} e^{\Psi + F}\omega_0^n.
\end{align*}
Note that $\Psi$ and $F$ are both bounded by $C(n, A)$, so we conclude that there exists a constant $C_p = C(n, A, \omega_0, p)>0$ such that 
\begin{equation}\label{eqn:new 10}
\int_X (\tr_{\omega_0}\omega)^{p-1 + \frac n{n-1}} \omega_0^n \le C_p \int_X (\tr_{\omega_0}\omega)^{p-1}\omega_0^n.
\end{equation}
When $p=2$ $$\int_X \tr_{\omega_0}\omega \omega_0^n = \int_X( n + \Delta_{\omega_0}\varphi )\omega_0^n = n \int_X \omega_0^n$$ is clearly bounded. Now we define a sequence $\{p_k\}$ with $p_0 = 2$ and $p_k = 2 + \frac{n}{n-1}k$. Then \eqref{eqn:new 10} implies that
$$\int_X (\tr_{\omega_0} \omega)^{p_k}\omega_0^n \le C_k \int_X (\tr_{\omega_0} \omega)^{p_{k-1}}\omega_0^n.$$
Since $p_k\to \infty$ as $k\to\infty$, iterating the inequality above gives that there exists a constant $C_k = C(n, A, \omega_0, k)>0$ such that 
$$\int_X (\tr_{\omega_0}\omega) ^{p_k}\omega_0^n\le C_k.$$
Lemma \ref{prop 1} then follows from this inequality and the H\"older inequality.

\begin{lemma}\label{cor 1}
For any $\beta\in (0,1)$, there exists a constant $C_\beta = C(n,A,\omega_0, \beta)>0$ such that 
$$\| \varphi\|_{C^{1,\beta}(X,\omega_0)}\le C_\beta.$$
\end{lemma}

{\it Proof.} By Lemma \ref{prop 1},  $f: = \Delta_{\omega_0} \varphi\in L^p(X,\omega_0^n)$ for any $p>0$. By the $W^{2,p}$-estimates for linear elliptic equations (c.f. Theorem 9.11 in \cite{GT}), we have
$$\| \varphi\|_{W^{2,p}(X,\omega_0)}\le C ( \| \varphi\|_{L^p(X,\omega_0^n)}  + \| f\|_{L^p(X,\omega_0^n)} )\le C_p.$$
The $C^{1,\beta}(X,\omega_0)$ bound of $\varphi$ then follows from the Sobolev embedding theorem (c.f. Corollary 7.11 in \cite{GT}) by taking $p>1$ sufficiently large. 

\begin{lemma}\label{lemma new 3}
The Ricci potential $u$ of $\omega = \omega_0+ i\ddbar \varphi$ satisfies 
$$\| u\|_{C^{\alpha}(X,\omega_0)}\le C_\alpha(n, A,\omega_0),$$ for any $\alpha\in (0,1)$.
\end{lemma}

{\it Proof.} Observe that 
$$|\nabla u|_{\omega_0}^2 \le |\nabla u|_{\omega}^2 \tr_{\omega_0} \omega \le A \tr_{\omega_0}\omega.$$
By Lemma \ref{prop 1}, it follows that $|\nabla u|_{\omega_0}\in W^{1,p}(X,\omega_0)$ for any $p>1$. The lemma then  follows from the Sobolev embedding theorem  by taking $p>1$ sufficiently large. 

\bigskip

To prove the $C^{2,\alpha}$-estimate of $\varphi$, we need the following recent result of Li-Li-Zhang, which weakens the condition of Y. Wang's result \cite{W} on the regularity assumption of $\varphi$. 

\begin{lemma}[\cite{LLZ} Theorem 1.2]
\label{LLZ lemma}
Let $B_2\subset \mathbb C^n$ be the Euclidean ball with radius $2$ and center $0$. Suppose $\varphi\in PSH(B_2)\cap C(B_2)$ solves the complex MA equation $$\det \varphi_{i\bar j} = f, \text{ in }B_2$$ with $f\ge \lambda >0$ for some positive $\lambda\in\mathbb R$ and $f\in C^\alpha(B_2)$ for some $\alpha\in (0,1)$. If $\varphi \in C^{1,\beta}(B_2)$ for some $\beta> 1- \frac{\alpha}{n(2+\alpha) - 1}$, then $\varphi\in C^{2,\alpha}(B_1)$ and the $C^{2,\alpha}(B_1)$-norm of $\varphi$ depends only on $n, \alpha, \beta,\lambda$, $\|\varphi\|_{C^{1,\beta}(B_2)}$ and $\| f\|_{C^\alpha(B_2)}$.
\end{lemma}

We arrive now at the $C^{2,\alpha}$ estimates for $\varphi$:

\begin{lemma}
\label{C2-alpha}
Under the conditions spelled out in the statement of Lemma \ref{lemma:C3}, there exists $\alpha>0$ with
\bea
\|\varphi\|_{C^{2,\alpha}}\leq C(n,A,\alpha)
\eea
for some constant $C(n,A,\o_0)$.
\end{lemma}

{\it Proof.} We note that by Lemma \ref{calpha} and Lemma \ref{lemma new 3}, the function on the right hand side of \eqref{eqn:MA} has uniform $C^{0,\alpha'}(X,\omega_0)$ estimate. Lemma \ref{cor 1} provides the $C^{1,\beta}(X,\omega_0)$ estimates of the K\"ahler potential $\varphi$. Then Lemma \ref{LLZ lemma} proves the $C^{2,\alpha}(X,\omega_0)$ estimates of $\varphi$. Q.E.D.

\bigskip

The following lemma is the key lemma that we shall need later for the proof of Theorem \ref{main} and Theorem \ref{main1}. It is
an immediate consequence of the $C^{2,\alpha}(X,\omega_0)$-estimates of $\varphi$, and the fact that the right hand side $e^F$ of the Monge-Amp\`ere equation (\ref{eqn:MA}) is bounded above and below:

\begin{lemma}\label{cor ricci potential}
There exists a constant $C=C(n, A, \omega_0)\ge 1$ such that 
$$ C^{-1}\omega_0\le \omega\le C \omega_0,$$ and $u\in C^{1,\alpha} (X,\omega_0)$ for any $\alpha\in (0,1)$.
\end{lemma}

\bigskip
Finally, we can complete the proof of Lemma \ref{lemma:C3}.
By Lemma \ref{C2-alpha}, the metric $g_{\bar ji}$ has uniform $C^{2,\alpha}$ norm. By Lemma \ref{lemma new 3} and Cramer's rule, its inverse $g^{i\bar j}$ also has uniform $C^{0,\alpha}(X,\omega_0)$ norm. The equation that $R_\o-n= \Delta_\omega u$ can be written locally in holomorphic coordinates as
$$g^{i\bar j}\frac{\partial ^2 u }{\partial z_i\partial \bar z_j}  = R_\o-n \in L^\infty.$$ Then the $C^{1,\alpha}(X,\omega_0)$-norm of $u$ follows from the $W^{2,p}$-estimates and Sobolev embedding theorem (c.f. \cite{GT}).

\medskip

Finally, once we have the $C^{1,\alpha}(X,\omega_0)$-norm of $u$, we can take $\frac{\partial}{\partial z_i}$ on both sides of  the equation \eqref{eqn:MA} and apply local Schauder estimates  to conclude that 
\begin{equation}\label{eqn:final}\| \varphi\|_{C^{3,\alpha}(X,\omega_0)} \le C(n, A,\omega_0,\alpha).\end{equation}
The proof of Lemma \ref{lemma:C3} is complete.

\section{Proof of Theorem 1}
\setcounter{equation}{0}

One direction in Theorem \ref{main} is a direct consequence of known results. Assume that $\lambda(X;A)>0$ for any $A$. 
By the work of Perelman (see \cite{ST} for a detailed account), for any given initial data in $c_1(X)$, the orbit of the K\"ahler-Ricci flow lies in a set $c_1(X;A)$ for some $A>0$. Thus a positive lower bound for $\lambda(X;A)$ implies a positive lower bound for the eigenvalue $\lambda(\om)$ along the K\"ahler-Ricci flow. By the results of \cite{PSSW1, ZZ}, the flow converges then to a K\"ahler-Einstein metric. 

\smallskip
The main issue in the present paper is to establish the other direction, namely that $\lambda(X;A)>0$ for any $A>0$ if a K\"ahler-Einstein metric is assumed to exist. But then the conditions 
of Lemma \ref{calpha} are satisfied, and for any fixed $A>0$, the metrics $\o\in c_1(X;A)$ have potentials which are uniformly bounded in $C^\alpha$-norm for some fixed $\alpha>0$. By Lemma \ref{cor ricci potential}, they are all equivalent. The desired bound for $\lambda(X;A)$ is then a consequence of the following lemma, which was essentially proved in \cite{PSSW0}, Lemma 1:

\begin{lemma}
\label{orthogonality}
Let $\o$, $\tilde\o$ be two metrics in $c_1(X)$ which are equivalent, in the sense that
\bea
\label{equivalent}
\kappa^{-1}\o\leq \tilde\o\leq \kappa\o
\eea
for some constant $\kappa>0$. Let $\lambda_\o$ and $\lambda_{\tilde \o}$ be the corresponding eigenvalues, as defined in (\ref{lambda}). Then 
\bea
c(\kappa,n)^{-1} \lambda_\o\leq \lambda_{\tilde\o}\leq c(\kappa,n)\lambda_\o
\eea
for some constant $c(\kappa,n)>0$ depending only on $\kappa$ and the dimension $n$.
\end{lemma}

{\it Proof.} Since this lemma is essential for our considerations and since its proof is short, we include the proof for the reader's convenience. In the definition (\ref{lambda}) for $\lambda_\o$ and $\lambda_{\tilde\o}$, the norms $\|\bar\p V\|_\o$ and $\|\bar\p V\|_{\tilde \o}$ as well as the norms $\|V\|_\o$ and $\|V\|_{\tilde \o}$ are already equivalent, since the metrics $\o$ and $\tilde\o$ are equivalent, and so are their volume forms $\o^n$ and $\tilde\o^n$. The main issue is the difference in the orthogonality conditions $\perp_\o$ and $\perp_{\tilde\o}$. 
To address this issue, consider any vector field $V$ with $V\perp_\o H^0(X,T^{1,0})$ and decompose it as
\bea
V=\tilde V+E
\eea
with $\tilde V\perp_{\tilde\o}H^0(X,T^{1,0})$ and $E\in H^0(X,T^{1,0})$. Taking inner products with respect to the metric $\o$ gives
\bea
0=\<\tilde V,E\>_\o+\<E,E\>_\o
\eea
and hence by the Cauchy-Schwarz inequality,
\bea
\|E\|_\o\leq \|\tilde V\|_\o.
\eea
We can now write for some constant $c_1(\kappa,n)$
\bea
\|\bar\p V\|^2_\o=\|\bar \p\tilde V\|^2_\o\geq c_1(\kappa,n)\|\bar\p \tilde V\|^2_{\tilde\o}
\eea
because $\o$ and $\tilde \o$ are equivalent, and at the same time, by the same equivalence and the triangle inequality,
\bea
\|V\|^2_\o\leq 2\|\tilde V\|^2_\o+2\|E\|^2_\o\leq 4\|\tilde V\|_\o\leq c_2(\kappa,n)\|\tilde V\|^2_{\tilde \o}.
\eea
It follows that
\bea
{\|\bar\p V\|_\o^2\over\|V\|_\o^2}
\geq 
{c_1(\kappa,n)\over c_2(\kappa,n)}
{\|\bar\p \tilde V\|_{\tilde\o}^2\over\|\tilde V\|_{\tilde \o}^2}
\geq
{c_1(\kappa,n)\over c_2(\kappa,n)}\lambda_{\tilde \o}
\eea
and hence $\lambda_\o\geq {c_1(\kappa,n)\over c_2(\kappa,n)}\lambda_{\tilde \o}$. Reversing the roles of $\o$ and $\tilde\o$ gives the inequality in the opposite direction. The lemma is proved, completing the proof of Theorem \ref{main}.

\section{Proof of Theorem 2}
\setcounter{equation}{0}

Again, one direction of the theorem follows from the results of Perelman and \cite{PSSW, ZZ}. To prove the other direction, namely that the existence of a K\"ahler-Einstein metric implies a strictly positive gap $\mu(X;A)$ for any $A>0$, we argue by contradiction. Recall the operator $L_\o$ defined for a metric $\o$ with Ricci potential $u$ by
$L_\o f\ = \ -g^{j\bar k}\nabla_j\nabla_{\bar k}f+g^{j\bar k}\nabla_{\bar k}f\nabla_ju \ - \ f$ and whose eigenvalues and eigenfunctions satisfy the identity
(\ref{futaki}).

Assume then that $X$ is K\"ahler-Einstein, and that there exists a sequence of metrics $\o_j=\o_0+i\p\bar\p \varphi_j\in c_1(X;A)$ such that the eigenvalues $\mu_j$ of the operator $L_{\omega_j}$ goes to $0$ as $j\to\infty$. We take $f_j$ to be eigenfunctions of $L_{\omega_j}$ with eigenvalues $\mu_j$, normalized by $\| f_j\|_{L^2( X, e^{-u_j} \omega_j^n )} = 1$. It follows from straightforward calculation that for any holomorphic vector field $V\in H^0(X, T^{1,0} X)$ \begin{equation}\label{eqn:orthogonal}
\int_X \langle \nabla_{\omega_j} f_j , V  \rangle_{\omega_j} e^{-u_j} \omega_j^n = 0.
\end{equation}
By Lemma \ref{cor ricci potential} and Lemma \ref{C2-alpha}, we can apply the elliptic estimates to $f_j$, which satisfies the linear equation
$$ - g_j^{p\bar q} \nabla_p\nabla_{\bar q} f_j + g_j^{p\bar q} \nabla_{\bar q} f_j \nabla_p u_j - f_j = \mu_j f_j  $$
to conclude that $$\| f_j \|_{C^{2,\alpha}(X, \omega_0)}\le C (n, A),\quad \forall\, j.$$
Up to a subsequence, we may assume the Ricci potentials $u_j$ converge in $C^{1,\alpha}$ to a function $u_\infty\in C^{1,\alpha}$, the metrics $\omega_j$ converge in $C^{1,\alpha}$ to a metric $\omega_\infty\in C^{1,\alpha}$, and the functions $f_j$ converge in $C^{2,\alpha}$ to a function $f_\infty \in C^{2,\alpha}$. In particular, we have $\| f_\infty\|_{L^2( X, e^{-u_\infty}\omega_\infty^n  )} = 1$. Passing to the limit, \eqref{eqn:orthogonal} gives that
\begin{equation}\label{eqn:orth1}\int_X \langle \nabla_{\omega_\infty} f_\infty, V \rangle_{\omega_\infty} e^{-u_\infty}\omega_\infty^n = 0,\quad \forall \, V\in H^0(X, T^{1,0} X) .\end{equation}

Observe that the equations (i.e. \eqref{futaki})
\begin{equation*}
\int_X |\bar \nabla \bar \nabla f_j|_{\omega_j}^2 e^{-u_j }\omega_j^n + \int_{X} |\bar \nabla f_j|^2_{\omega_j}  e^{-u_j}\omega_j^n = (1+\mu_j) \int_X |\bar \nabla f_j|^2_{\omega_j}  e^{-u_j}\omega_j^n
\end{equation*}
hold for any $j$. Since  $\mu_j\to 0$, passing to limit we get
$$\int_X |\bar \nabla \bar \nabla f_\infty|^2_{\omega_\infty} e^{-u_\infty}\omega_\infty^n = 0,$$ which implies $\nabla \nabla f_\infty = 0$, i.e. $\nabla_{\omega_\infty} f_\infty$ is a holomorphic vector field. From \eqref{eqn:orth1} we conclude that $\int_X |\nabla f_\infty|^2_{\omega_\infty} e^{-u_\infty} \omega_\infty^n = 0$. However, this contradicts the identity 
$$1= \int_X f_\infty^2 e^{-u_\infty} \omega_\infty^n = \int_X |\nabla f_\infty|^2_{\omega_\infty} e^{-u_\infty}\omega_\infty^n.$$
The proof of Theorem \ref{main1} is complete.

\medskip
We observe that this argument could have been used also for the proof of Theorem \ref{main}. However, the argument there is more direct, and provides more precise information on the bounds for $\lambda_\o$.

\section{Further remarks}
\setcounter{equation}{0}

We note that in Theorem \ref{main}, we cannot in general replace the gap $\lambda(X;A)$ for each $A>0$ by the gap $\lambda(X)={\rm inf}_{\om\in c_1(X)}\lambda_\om$ over all of $c_1(X)$. A simple counterexample is provided by the $2$-dimensional sphere, which admits a K\"ahler-Einstein metric, but can be seen to have
\bea
\lambda(S^2)=0
\eea
as follows. Let $\eta:{\bf R}\to{\bf R}$ be a smooth increasing function such that 
$\eta=0$ on $(-\i,1/3]$ and $\eta=1$ on $[2/3,\i)$. Let $a,N>0$ and let $f: [0,3N+2]\to {\bf R}$ be a non-negative concave function, positive and smooth on $(0,3N+2)$ such that 
\begin{enumerate}
\item $f(0)=f(3N+2)=0$
\item $f(x)=a$ for $x\in [1,3N+1]$
\end{enumerate}
and let $X$ be the surface obtained by revolving the graph of $y=f(x)$ around the $x$ axis. Thus $X$ is a smooth manifold (if we choose $f$ so that its tangent line is vertical at $0$ and $3N+2$ and is tangent to the graph to infinite order), looks like a cigar, is flat between $x=1$ and $x=3N+1$ and is diffeomorphic to $S^2$. Moreover, we can choose $a$ so that the area of $X$ is $1$ (so $a$ is roughly $1\over 3N\cdot 2\pi$). Let $g_N$ be metric obtained by restricting the euclidean metric in ${\bf R}^3$. Let $V_1$ be a smooth vector field on $X$
defined as follows.
$$ 
V_1= \eta(x-1)\eta(N+1-x){\p\over\p x}
$$
so $V_1$ is a smooth vector field on $X$ compactly supported in 
$\{(x,y,z)\in M\,:\, x\in (1,N+1)\}$. Similarly we define $V_2$ supported in 
$(N+1,2N+1)$ and $V_3$ supported in $(2N+1,3N+1)$.

Next we let 
\be 
V= c_1V_1+c_2V_2+c_3V_3
\ee
where the $c_i\in{\bf R}$ are chosen so that $V$ is orthogonal to the 3-dimensional space of holomorphic vector fields and $c_1^2+c_2^2+c_3^2=1$. Now $|V_i|$ is roughly equal to $1$ so $\|V_i\|_{L^2}\sim 1/3$ so 
$\|V\|_{L^2}\sim c_1^2\|V_1\|^2_{L^2}+c_2^2\|V_2\|^2_{L^2}+c_3^2\|V_3\|^2_{L^2}\sim {1\over 9}$. On the other hand $\nabla V_1=0$ for $2<x<N$ so
$\|\nabla V_1\|_{L^2}^2 =O({1\over N})$ which implies $\|\nabla V\|_{L^2}^2 =O({1\over N})$. In particular, $\lambda_{\o_N} \leq O({1\over N})$. This establishes our claim.

\bigskip

\noindent Bin Guo, Department of Mathematics and Computer Sciences, Rutgers University, Newark, NJ  07102

Email: bguo@rutgers.edu

\medskip

\noindent Duong H. Phong, Department of Mathematics, Columbia University, New York, NY 10027

Email: phong@math.columbia.edu

\medskip

\noindent Jacob Sturm,  Department of Mathematics and Computer Sciences, Rutgers University, Newark, NJ  07102

Email: sturm@andromeda.rutgers.edu

\end{document}